\pdfoutput=1
\RequirePackage{ifpdf}
\ifpdf 
\documentclass[pdftex]{sigma}
\else
\documentclass{sigma}
\fi

\newtheorem{Theorem}{Theorem}[section]
\newtheorem{thm}{Theorem}[section]

\newtheorem{prop}[Theorem]{Proposition}
{\theoremstyle{definition}

}
\DeclareMathOperator{\Hom}{Hom}
\DeclareMathOperator{\Rep}{Rep}
\DeclareMathOperator{\Imop}{Im}
\begin{document}
\newcommand{\blambda}{{\mbox{\boldmath$\lambda$}}}
\newcommand{\bmu}{{\mbox{\boldmath$\mu$}}}
\newcommand{\bnu}{{\mbox{\boldmath$\nu$}}}

\newcommand{\arXivNumber}{1403.4750}

\allowdisplaybreaks

\renewcommand{\thefootnote}{$\star$}

\renewcommand{\PaperNumber}{058}

\FirstPageHeading

\ShortArticleName{Schur positivity and Kirillov--Reshetikhin modules}

\ArticleName{Schur Positivity and Kirillov--Reshetikhin Modules\footnote{This paper is a~contribution to the Special
Issue on New Directions in Lie Theory.
The full collection is available at
\href{http://www.emis.de/journals/SIGMA/LieTheory2014.html}{http://www.emis.de/journals/SIGMA/LieTheory2014.html}}}

\Author{Ghislain FOURIER~$^\dag$ and David HERNANDEZ~$^\ddag$}

\AuthorNameForHeading{G.~Fourier and D.~Hernandez}

\Address{$^\dag$~School of Mathematics and Statistics, University of Glasgow, UK}
\EmailD{\href{mailto:ghislain.fourier@glasgow.ac.uk}{ghislain.fourier@glasgow.ac.uk}}

\Address{$^\ddag$~Sorbonne Paris Cit\'e, Univ Paris Diderot-Paris 7, Institut de Math\'ematiques de Jussieu~-\\
\hphantom{$^\ddag$}~Paris Rive Gauche CNRS UMR 7586, B\^at. Sophie Germain, Case 7012,75205 Paris, France}
\EmailD{\href{mailto:hernandez@math.jussieu.fr}{hernandez@math.jussieu.fr}}

\ArticleDates{Received April 04, 2014, in f\/inal form May 29, 2014; Published online June 04, 2014}

\Abstract{In this note, inspired by the proof of the Kirillov--Reshetikhin conjecture, we consider tensor products of
Kirillov--Reshetikhin modules of a~f\/ixed node and various level.
We f\/ix a~positive integer and attach to each of its partitions such a~tensor product.
We show that there exists an embedding of the tensor products, with respect to the classical structure, along with the
reverse dominance relation on the set of partitions.}


\Classification{17B10; 17B37; 05E05}

\renewcommand{\thefootnote}{\arabic{footnote}} \setcounter{footnote}{0}

\section{Introduction}

This note is inspired by two results on certain modules of simple, f\/inite-dimensional complex Lie algebras.
The f\/irst one is an immediate consequence of the famous Clebsch--Gordan formula on decompositions of tensor product of
simple $\mathfrak{sl}_2(\mathbb{C})$-modules.
Namely let $V(m) = \operatorname{Sym}^m \mathbb{C}^2$, be the~$n$-th symmetric power of the natural representation, then
\begin{gather*}
V(n) \otimes V(m) \cong_{\mathfrak{sl}_2} V(n +m) \oplus V(n+m- 2) \oplus \dots \oplus V(n+m - 2\min\{n,m\}).
\end{gather*}
Which implies that for $m_1, m_2 \leq m$ there exists a~surjective map of $\mathfrak{sl}_2$-modules
\begin{gather*}
V(m_2) \otimes V(m- m_2) \twoheadrightarrow V(m_1) \otimes V(m - m_1)
\end{gather*}
if and only if $\min\{m_1, m - m_1 \} \leq \min \{m_2, m- m_2 \}$.
Using this inequality, we obtain an order~$\preceq$ on partitions of~$n$ of length~$2$, $\mathcal{P}(m,2)$.
By taking the point of view from symmetric functions, $s_m$~being the character of~$V(m)$, we have
\begin{gather*}
(m_1, m- m_1) \preceq (m_2, m - m_2) \Leftrightarrow s_{m_2} s_{m - m_2} - s_{m_1} s_{m - m_1} \in \sum\limits_{k \geq
0} \mathbb{Z}_{\geq 0} s_k.
\end{gather*}
As the characters are also known as Schur functions, this property of the left hand side is also known
as \textit{Schur positivity}.
A~generalization of this order to $\mathfrak{sl}_n(\mathbb{C})$ and further to a~simple f\/inite-dimensional Lie algebra
$\mathfrak{g}$ of arbitrary type was investigated in~\cite{DP07} (resp.~\cite{CFS13, Fou14}).

The other inspiration comes from certain character identities for classical limits of Kirillov--Reshetikhin modules for
$U_q(\hat{\mathfrak{g}})$, the untwisted quantum af\/f\/ine algebra associated to $\mathfrak{g}$, namely the $Q$-systems.
Kirillov--Reshetikhin modules $W_{m,a}^{(i)}$ are indexed by a~node of the Dynkin diagram, say $i \in I$, a~positive
level~$m$ and a~parameter $a \in \mathbb{C}^*$.
For more details on Kirillov--Reshetikhin modules and their importance in the category of f\/inite-dimensional
$U_q(\hat{\mathfrak{g}})$-modules we refer to~\cite{CH10}.

We denote by ${\rm KR}(m \omega_i)$ the $\mathfrak g$-module obtained through the limit $q \rightarrow 1$ of $W_{m,a}^{(i)}$,
note that the classical structure is independent of~$a$.
Further denote $\operatorname{char}_{\mathfrak g} {\rm KR}(m \omega_i)$ the classical character of ${\rm KR}(m \omega_i)$, then the
$Q$-system is the following character identity~\cite{hcr, Nad}
\begin{gather*}
\operatorname{char}_{\mathfrak g} {\rm KR}(m \omega_i) \operatorname{char}_{\mathfrak g} {\rm KR}(m \omega_i) =
\operatorname{char}_{\mathfrak g} {\rm KR}((m+1) \omega_i)\operatorname{char}_{\mathfrak g} {\rm KR}((m-1) \omega_i) +
\operatorname{char}_{\mathfrak g} S_{m}^{(i)},
\end{gather*}
where $ S_{m}^{(i)}$ denotes the classical limit of a~tensor product of certain Kirillov--Reshetikhin modules (depending
on~$i$ and~$m$).

From the $Q$-system one has the immediate consequence, that there exists a~surjective map of $\mathfrak g$-modules
\begin{gather*}
{\rm KR}(m \omega_i) \otimes {\rm KR}(m \omega_i) \twoheadrightarrow {\rm KR}((m+1) \omega_i) \otimes {\rm KR}((m-1) \omega_i).
\end{gather*}
Considering the partial order on partitions, we have $(m-1, m+1) \prec (m,m) \in \mathcal{P}(2m,2)$.
So for the maximal element and its predecessor this might be seen as a~generalization of Schur positivity (to Lie
algebras of arbitrary type).
Note that in~\cite{CV} it is proved that there exists a~surjective map of modules for the current algebra of $\mathfrak
g$, namely of fusion products of modules of Kirillov--Reshetikhin modules for the partitions $(m-1, m+1)$, $(m,m)$.
Their work is also motivated by the $Q$-system.

Combining the partial order on $\mathcal{P}(m,2)$ and this consequence of the $Q$-system relation was the starting point
of this paper.
We have generalized the arguments in the proof of the Kirillov--Reshetikhin conjecture, e.g.~the character of the Kirillov--Reshetikhin modules satisfy the $Q$-system.
Using this, we have proved that for all $\mathfrak{g}$ of arbitrary type and $(m_1, m- m_1) \preceq (m_2, m - m_2) \in
\mathcal{P}(m,2)$ there exists a~surjective map of $\mathfrak g$-modules
\begin{gather*}
{\rm KR}(m_2 \omega_i) \otimes {\rm KR}((m - m_2) \omega_i) \twoheadrightarrow {\rm KR}(m_1\omega_i)\otimes {\rm KR}((m - m_1) \omega_i).
\end{gather*}
This might be seen as a~generalization of the $\mathfrak{sl}_2$-case as well as of the consequence of the $Q$-system
property.

More generally, we associate to each partition $(m_1 \geq m_2 \geq \dots \geq m_k > 0)$ of~$m$ a~tensor product of
Kirillov--Reshetikhin modules
\begin{gather*}
{\rm KR}(m_1 \omega_i) \otimes \dots \otimes {\rm KR}(m_k \omega_i).
\end{gather*}
By considering the reverse dominance relation on partitions of~$m$, $\mathcal{P}(m)$, we can show further
(using that 
the cover relation is induced by the cover relations on the set of partitions of length~$2$)
that if $(m_1 \geq \dots
\geq m_{k_1} > 0) \preceq (n_1 \geq \dots \geq n_{k_2} > 0)$, then there exists a~surjective map of
$\mathfrak{g}$-modules
\begin{gather*}
{\rm KR}(n_1 \omega_i) \otimes \dots \otimes {\rm KR}(n_{k_2} \omega_i) \twoheadrightarrow {\rm KR}(m_1 \omega_i) \otimes \dots \otimes
{\rm KR}(m_{k_1} \omega_i).
\end{gather*}
In the $\mathfrak{sl}_n$-case, the last statement can be deduced from the results on row shuf\/f\/les in~\cite{LPP07}.
Even more general, in the case where $\omega_i$ is minuscule, e.g.~${\rm KR}(m \omega_i)$ is a~simple $\mathfrak{g}$-module,
this last statement has been proved in~\cite{CFS13}.
The authors were constructing an explicit bijection of the highest weight vectors in terms of LS-paths.
Our approach avoids these combinatorics.

Since ${\rm KR}(m \omega_i)$ can be also constructed as a~module for $\mathfrak g \otimes \mathbb{C}[t]$ (by using some
``evaluation parameter'' $a \in \mathbb{C}$), see Section~\ref{chari}, one might ask if there is a~surjection also as
$\mathfrak g \otimes \mathbb{C}[t]$-modules.
The natural object to be considered here is the fusion product introduced in~\cite{FL99}.
This is the associated graded module (with respect to the degree f\/iltration on $U(\mathfrak g \otimes \mathbb{C}[t])$)
of the tensor product of the Kirillov--Reshetikhin modules with pairwise distinct evaluations.
Can the surjection in Theorem~\ref{two-factor} be actually obtained from a~surjective map of the corresponding fusion products?

We should remark here that a~similar result on Schur positivity on tensor products of simple $\mathfrak{g}$-modules of
arbitrary highest weight~$\lambda$ was conjectured in~\cite{DP07} and~\cite{CFS13} (see Section~3.5 for
more details).
Our result suggests, that this generalized Schur positivity may hold along the partial order on tensor products of $q
\mapsto 1$ limits of minimal af\/f\/inizations of~$V(\lambda)$ (the ``minimal'' module of the quantum af\/f\/ine algebra having
a~simple quotient whose limit is isomorphic to $V(\lambda)$, see~\cite{CP96} for more details).

\looseness=-1
In the $\mathfrak{sl}_n$-case, the restriction to $\mathfrak{sl}_n$ of the limit of such a~minimal af\/f\/inization is
nothing but the simple $\mathfrak{sl}_n$-module~$V(\lambda)$, so this is the conjecture of Schur positivity by Lam,
Postnikov and Pylyavskyy, cited in~\cite{DP07}.
For other types, the limit of a~minimal af\/f\/inization is not a~simple $\mathfrak g$-module in general, for example
Kirillov--Reshetikhin modules are minimal af\/f\/inizations of~$V(m \omega_i)$.
It might be interesting to investigate on minimal af\/f\/inizations of other than rectangular weights.

In Section~\ref{section-pre} we brief\/ly recall the reverse dominance relation on partitions and some basics on
Kirillov--Reshetikhin modules.
In Section~\ref{section-main} we state the main theorem, while the proof follows in Section~\ref{section-proof}.

\section{Preliminaries}
\label{section-pre}
Let $\mathfrak g$ be a~f\/inite-dimensional simple, complex Lie algebra of rank~$n$ and Cartan matrix~$C$.
Let $\mathfrak n^+ \oplus \mathfrak h \oplus \mathfrak n^-$ be a~triangular decomposition.
We denote the set of (positive) roots~$R$ ($R^+$ resp), the (dominant) integral weights~$P$ (resp.~$P^+$).
We denote the simple roots $\{\alpha_1, \dots, \alpha_n \}$, the fundamental weight $\{\omega_1, \dots, \omega_n\}$, $I
= \{1, \dots, n \}$.
For every $\alpha \in R^+$ we f\/ix a~$\mathfrak{sl}_2$-triple $\{e_{\alpha}, h_\alpha, f_{\alpha} \}$.

\subsection{Partial order}
We recall the reverse dominance order on partitions.
For this, let $m \geq 1$ be positive integer and denote by $\mathcal{P}(m)$ the set of partitions of~$m$:
\begin{gather*}
\mathcal{P}(m) = \{(m_1 \geq \dots \geq m_k > 0)  \mid  m_j \in \mathbb{Z}
~
\text{and}
~
m_1 + \dots + m_k = m  \}.
\end{gather*}
This is a~f\/inite set and the reverse dominance relation on $\mathcal{P}(m)$ is def\/ined as follows:
\\
Let $\blambda = (m_1 \geq \dots \geq m_{k_1} > 0)$, $\bmu = (n_1 \geq \dots \geq n_{k_2} > 0) \in \mathcal{P}(m)$.
Then
\begin{gather*}
\blambda \preceq \bmu
\quad
:\Leftrightarrow
\quad
\forall\,  j = 1, \dots, \min\{k_1, k_2 \}:
\
m_1 + \dots + m_j \geq n_1 + \dots + n_j.
\end{gather*}
Obviously, \looseness=-1 this gives a~partial order on $\mathcal{P}(m)$ with a~smallest element $(m > 0)$ and largest element $(1 \geq
1 \geq \dots \geq 1 > 0)$.
Moreover, if we consider partitions of a~f\/ixed length~$k$ only, $\mathcal{P}(m, k)$, then there is also a~unique maximal
element.
Namely if $m = \ell k + p$, where $0 \leq p < k$, then
\begin{gather*}
\blambda = (\ell + 1 \geq \dots \geq \ell + 1 > \ell \geq \dots \geq \ell)
\end{gather*}
is the unique maximal element in $\mathcal{P}(m,k)$.

\subsection{Cover relation}
We recall the notion of the cover relation induced by $\preceq$, e.g.~we say $\bmu$ covers $\blambda$ if 
\begin{enumerate}\itemsep=0pt

\item[1)] $\blambda \preceq \bmu$ and

\item[2)] $\blambda \preceq \bnu \preceq \bmu$ implies $\bnu= \blambda$ or $\bnu = \bmu$.
\end{enumerate} 
Since $\mathcal{P}(m)$ is a~f\/inite set, we can f\/ind for each pair $\blambda \preceq \bmu$ partitions $\bnu_0, \dots,
\bnu_{\ell}$ such that
\begin{gather*}
\blambda = \bnu_0 \preceq \bnu_1 \preceq \dots \preceq \bnu_{\ell} = \bmu
\end{gather*}
and $\bnu_{i}$ covers $\bnu_{i-1}$ for all~$i$.
To understand the partial order on $\mathcal{P}(m)$ it is therefore suf\/f\/icient to understand the cover relation on
$\mathcal{P}(m)$.
The following proposition was proved in~\cite[Proposi\-tion~3.5]{CFS13} (for simplicity of notation we assume that
$\blambda$, $\bmu$ have the same length by adding~$0$~parts to at most one of the both).

\begin{prop}
\label{cover}
Let $\blambda = (m_1 \geq \dots \geq m_{k} \geq 0)$, $\bmu = (n_1 \geq \dots \geq n_{k} \geq 0) \in \mathcal{P}(m)$ and
$n_k$ or $m_k \neq 0$.
Suppose $\bmu$ covers $\blambda$.
Then there exists $i < j$ such that
\begin{gather*}
n_{\ell} =
\begin{cases}
m_\ell & \text{if}\quad  \ell \neq i, j,
\\
m_{\ell} - 1 & \text{if}\quad  \ell = i,
\\
m_{\ell} +1 & \text{if}\quad  \ell = j.
\end{cases}
\end{gather*}
The cover relation on partitions $\mathcal{P}(m)$ is completely determined by the cover relation on partitions of length
$2$.
\end{prop}

\subsection{Quantum af\/f\/ine algebras and their representations}

We give a~brief reminder on quantum af\/f\/ine algebras and their f\/inite-dimensional representations.
For more details we refer to~\cite{CH10, hcr}.

Let $q\in\mathbb{C}^*$ which is not a~root of unity.
Let $U_q(\hat{\mathfrak{g}})$ be the untwisted quantum af\/f\/ine algebra associated to $\mathfrak{g}$.
The simple objects of the category $\mathcal{C}$ of (type~$1$) f\/inite-dimensional representations of
$U_q(\hat{\mathfrak{g}})$ are parametrized by dominant monomials of the ring $\mathcal{Y} = \mathbb{Z}[Y_{i,a}^{\pm
1}]_{1\leq i\leq n,a\in\mathbb{C}^*}$, that is for each such monomial $m = \prod\limits_{i\in I,
a\in\mathbb{C}^*}Y_{i,a}^{u_{i,a}}$ which is dominant (the $u_{i,a}\geq 0$), there is a~corresponding simple object
$L(m)$ in $\mathcal{C}$.
For example, for $i\in I$, $m\geq 0$, $a\in\mathbb{C}^*$, we have the Kirillov--Reshetikhin module
\begin{gather*}
W_{m,a}^{(i)} = L\big(Y_{i,a}Y_{i,aq_i^2}\cdots Y_{i,aq_i^{2(m-1)}}\big).
\end{gather*}
Here $q_i = q^{r_i}$ where $r_i$ is the length of simple root $\alpha_i$.

The~$q$-character morphism~\cite{Fre} is an injective ring morphism def\/ined on the Grothendieck ring
$\Rep(U_q(\hat{\mathfrak{g}}))$ of the tensor category $\mathcal{C}$:
\begin{gather*}
\chi_q: \ \Rep(U_q(\hat{\mathfrak{g}}))\rightarrow \mathcal{Y}.
\end{gather*}
\begin{thm}[\protect{\cite{Fre2, Fre}}] For~$m$ a~dominant monomial, we have
\begin{gather*}
\chi_q(L(m))\in m\mathbb{Z}_{\geq 0}\big[A_{i,a}^{-1}\big]_{1\leq i\leq n, a\in\mathbb{C}^*},
\end{gather*}
where
\begin{gather*}
A_{i,a}=Y_{i,aq_i^{-1}}Y_{i,aq_i}\prod\limits_{\{j|C_{j,i}=-1\}}Y_{j,a}^{-1}
\prod\limits_{\{j|C_{j,i}=-2\}}Y_{j,aq^{-1}}^{-1}Y_{j,aq}^{-1}
\prod\limits_{\{j|C_{j,i}=-3\}}Y_{j,aq^2}^{-1}Y_{j,a}^{-1}Y_{j,aq^{-2}}^{-1}.
\end{gather*}
An element in $\Imop(\chi_q)$ is characterized by the multiplicity of its dominant monomials.
\end{thm}

Note that we have a~partial ordering on the monomials of $\mathcal{Y}$: $m\preceq m'$ if $m'm^{-1}$ is a~product of
monomials $A_{i,a}$.
The f\/irst statement in the theorem can be reformulated by saying that all monomials in $\chi_q(L(m))$ are lower than~$m$ 
for this ordering, that is~$m$ is the highest monomial.

As consequence of the second statement, if we know that the~$q$-character of a~simple module has a~unique dominant
monomial, its~$q$-character can be reconstructed (this is the Frenkel--Mukhin algorithm~\cite{Fre2}).
This property has been proved in the important case of Kirillov--Reshetikhin modules, which led to the proof of the
Kirillov--Reshetikhin conjecture.
It was f\/irst proved by Nakajima~\cite{Nad} for $ADE$-types and in~\cite{hcr} with a~dif\/ferent proof which can be
extended to the general case.

\begin{thm}[\protect{\cite{hcr, Nad}}]
\label{tsys}
The~$q$-character of Kirillov--Reshetikhin module has a~unique dominant monomial.
This implies the~$T$-system in the Grothendieck ring $\Rep(\mathcal{U}_q(\hat{\mathfrak{g}}))$:
\begin{gather*}
\big[W_{m,a}^{(i)}\otimes W_{m,aq_i^2}^{(i)}\big] = \big[W_{m+1,a}^{(i)}\otimes W_{m-1,aq_i^2}^{(i)}\big] + \big[S_{m,a}^{(i)}\big],
\end{gather*}
where $S_{m,a}^{(i)}$ is a~tensor product of Kirillov--Reshetikhin modules.
The representations $S_{m,a}^{(i)}$ and $W_{m+1,a}^{(i)}\otimes W_{m-1,aq_i^2}^{(i)}$ are simple.
\end{thm}

\subsection{The classical limit}
Let ${\rm KR}(m\omega_i)$ be the restriction of $W_{m,a}^{(i)}$ as a~$U_q(\mathfrak{g})$-module (it is well known it does not
depend on~$a$, see references in~\cite{hcr}).
We denote by the same symbol its limit at $q = 1$ (that is we consider the corresponding $\mathfrak{g}$-module).

${\rm KR}(m\omega_i)$ decomposes into a~direct sum of f\/initely many simple $\mathfrak{g}$-modules.
This decomposition is computed for $\mathfrak{g}$ of classical type, namely type~$A$, $B$, $C$, $D$ as well as for certain nodes
for exceptional types (see \cite[Theorem~1]{Cha01} and~\cite[Chapter~3]{Cha01}).
For more on decompositions of Kirillov--Reshetikhin modules see~\cite{HKOTT}.

\subsection{Chari modules}
\label{chari}
For the readers convenience we should clarify the relation to \textit{Kirillov--Reshetikhin modules for current
algebras}, although it is not used in this note at all:

Denote $\mathfrak g \otimes \mathbb{C}[t]$ the current algebra of $\mathfrak g$.
In~\cite{Cha01}, f\/inite-dimensional modules for $\mathfrak g \otimes \mathbb{C}[t]$ were introduced.
These are supposed to be classical analogs of Kirillov--Reshetikhin modules.
Namely for $m \in \mathbb{Z}_{\geq 0}$, $a\in \mathbb{C}$, $i \in I$, $C(m \omega_i,a)$ is the $\mathfrak g \otimes
\mathbb{C}[t]$-module generated by a~non-zero vector~$w$ subject to the relations
\begin{gather*}
\mathfrak n^+ \otimes \mathbb{C}[t].w = 0,
\qquad
(h \otimes 1).w = m \omega_i(h) w,
\qquad
\mathfrak h \otimes (t-a) \mathbb{C}[t].w = 0,
\\
(f_{\alpha})^{m \omega_i(h_\alpha) + 1}.w = 0,
\qquad
(f_{\alpha_i} \otimes (t-a)).w = 0.
\end{gather*}

The following was proven for $\mathfrak g$ of classical type, namely of type $A$, $B$, $C$, $D$ in~\cite[Corollary~2.1]{Cha01}.
The general case (e.g., an arbitrary simple, f\/inite-dimensional complex Lie algebra) can be deduced as the special case
of a~single tensor factor of~\cite[Corollary~5.1]{Ked11}, this has been proved by using a~pentagram of identities
see~\cite[Section~1.2]{Ked11}.

\begin{thm}
For $i \in I$ and $m \geq 0$ we have
\begin{gather*}
C(m \omega_i, a) \cong_{\mathfrak g} {\rm KR}(m \omega_i).
\end{gather*}
\end{thm}

\section{Main result}
\label{section-main}

{\bf 3.1.~Main theorem.}
For $\blambda = (m_1 \geq \dots \geq m_k > 0) \in \mathcal{P}(m)$ and $i \in I$ we denote the tensor product of the
associated KR modules as:
\begin{gather*}
{\rm KR}(\blambda, i):= {\rm KR}(m_1\omega_i) \otimes \dots \otimes {\rm KR}(m_k\omega_i).
\end{gather*}
Note that ${\rm KR}(\bf{0}, \omega_i) \cong \mathbb{C}$ is the trivial representation.
We are mainly interested in the $\mathfrak{g}$-structure of these modules.
Since this is the tensor product of f\/inite-dimensional $\mathfrak{g}$-modules, it decomposes into the direct sum of
simple modules.
Due to the structure of the Kirillov--Reshetikhin modules, we see immediately that the maximal weight occurring is $m
\omega_i$ and that the corresponding weight space has dimension~$1$.

The main purpose is of this paper is the study of the $\mathfrak{g}$-structure of these modules along the partial order
and the main theorem is the following
\begin{thm}
\label{main-thm}
Let $m \geq 1$, $i \in I$ and $\blambda \preceq \bmu \in \mathcal{P}(m)$, then
\begin{gather*}
\dim (\Hom_{\mathfrak g}({\rm KR}(\blambda, i), V(\tau))) \leq \dim (\Hom_{\mathfrak g}({\rm KR}(\bmu, i), V(\tau)))
\end{gather*}
for all $\tau \in P^+$.
\end{thm}
In other words, there exists a~surjective map of $\mathfrak g$-modules
\begin{gather*}
{\rm KR}(\bmu, i) \longrightarrow {\rm KR}(\blambda, i).
\end{gather*}
\begin{proof}
Let $\blambda \preceq \bmu  \in \mathcal{P}(m)$.
We denote $\bmu = (n_1 \geq \dots \geq n_{k_2} \geq 0)$.
It is suf\/f\/icient to prove the statement for the case where $\bmu$ covers $\blambda$.
In this case, we know by Proposition~\ref{cover} there exists $p < q$ such that $n_{\ell} = m_{\ell}$ for all $\ell \neq
p,q$.
This implies that
\begin{gather*}
{\rm KR}(\bmu, i) = {\rm KR}((m_p - 1)\omega_i) \otimes {\rm KR}((m_q +1) \omega_i) \otimes \bigotimes_{\ell \neq p,q} {\rm KR}(m_{\ell}
\omega_i).
\end{gather*}
So to prove the statement it is enough to give a~proof for partitions of~$m$ of length~$2$.
So let $ m_1 \geq m_2 > 0$, we have to show that we have a~surjective map of $\mathfrak g$-modules
\begin{gather*}
{\rm KR}(m_1 \omega_i) \otimes {\rm KR}(m_2 \omega_i) \twoheadrightarrow {\rm KR}((m_1 +1) \omega_i) \otimes {\rm KR}((m_2 - 1) \omega_i).
\end{gather*}
This will be proven in the next section, Theorem~\ref{two-factor}.
\end{proof}

{\bf 3.2.}~Before proving the last ingredient in the proof (Theorem~\ref{two-factor}), we shall make a~couple of remarks.
First of all, Theorem~\ref{main-thm} was proved in~\cite[Theorem ii)]{CFS13} for the special case where $\omega_i$ is
a~minuscule weight.
In this case, ${\rm KR}(m \omega_i) \cong V(m \omega_i)$ as a~$\mathfrak g$-module (see, e.g.,~\cite{Cha01}).
This case covers the $\mathfrak{sl}_{n+1}$ case as well as certain special cases for other types.
The proof given there uses the combinatorics of LS paths.
Namely, an injection on the level of paths in the tensor product is given.
This proof does not extend to other cases since the combinatorics of LS-paths are more complicated for non-minuscule
weights.

{\bf 3.3.}~It would be interesting to study the kernel of the map ($\blambda \preceq \bmu \in \mathcal{P}(m)$)
\begin{gather*}
{\rm KR}(\bmu, i) \longrightarrow {\rm KR}(\blambda, i).
\end{gather*}
In the formerly studied case of $Q$-systems, this kernel turns out to be a~tensor product of Kirillov--Reshetikhin modules
again (Theorem~\ref{tsys}).
This is not true in general: for $\mathfrak{sl}_n$, $i = 2$ and the quotient
\begin{gather*}
{\rm KR}(5,2) \otimes {\rm KR}(1,2) \longrightarrow {\rm KR}(6,2),
\end{gather*}
the kernel is not a~(non-trivial) tensor product at all.
It is also not a~minimal af\/f\/inization as they are all simple $\mathfrak{sl}_n$-modules.

{\bf 3.4.}~The partial order on $\mathcal{P}(m)$ is a~special case of the more general poset $\mathcal{P}(\lambda)$ introduced
in~\cite{CFS13}.
Here $\lambda \in P^+$ is a~dominant weight and the elements in the set are partitions of~$\lambda$, namely $\blambda =
(\lambda_1, \dots, \lambda_k \geq 0)$, with $\lambda_i \in P^+$ and $\lambda_1 + \dots + \lambda_k = \lambda$.
Then the partial order on $\mathcal{P}(\lambda)$ is given by: $\blambda \preceq \bmu $ if and only if for all $\alpha
\in R^+$ and $\ell \geq 1$:
\begin{gather*}
\min_{i_1 < \dots < i_\ell} \{(\lambda_{i_1} + \dots + \lambda_{i_\ell})(h_\alpha)\} \leq \min_{i_1 < \dots < i_\ell}
\{(\mu_{i_1} + \dots + \mu_{i_\ell})(h_\alpha)\}.
\end{gather*}
Then~\cite[Theorem 1(i)] {CFS13} gives for $\blambda \preceq \bmu$
\begin{gather*}
\dim (V(\lambda_1) \otimes \dots \otimes V(\lambda_{k_1})) \leq \dim (V(\mu_1) \otimes \dots \otimes V(\mu_{k_2})).
\end{gather*}

{\bf 3.5.}~It is conjectured~\cite[Conjecture~1]{CFS13} that if $(\lambda_1, \lambda_2) \preceq (\mu_1, \mu_2) \in
\mathcal{P}(\lambda, 2)$ then
\begin{gather*}
\dim (\Hom_{\mathfrak g} (V(\lambda_1) \otimes V(\lambda_2), V(\tau))) \leq \dim (\Hom_{\mathfrak g}(V(\mu_1) \otimes
V(\mu_2), V(\tau)))
\end{gather*}
for all $\tau \in P^+$.

This conjecture was made before for $\mathfrak{sl}_{n+1}$ by Lam, Postnikov, Pylyavskyy (cited in~\cite[Conjec\-ture~1]{DP07}).
In this case, the conjecture is equivalent to the following statement on the level of characters:
\begin{gather*}
s_{\mu_1} s_{\mu_2} - s_{\lambda_1} s_{\lambda_2} = \sum\limits_{\tau \in P^+} c_{\tau} s_{\tau},
\qquad
\text{and}
\qquad
c_\tau \geq0
\quad
\forall
\,
\tau \in P^+,
\end{gather*}
where $s_\lambda$ is the Schur function corresponding to~$\lambda$, e.g.
the character of the simple $\mathfrak{sl}_{n+1}$-mo\-du\-le~$V(\lambda)$.
In other words, it is conjectured that the dif\/ference of the products of the corresponding Schur functions is
\textit{Schur positive}.

A big step forward in proving this conjecture has been made with~\cite[Theorem~1]{DP07}:
\begin{gather*}
\dim (\Hom_{\mathfrak g} (V(\lambda_1) \otimes V(\lambda_2), V(\tau))) \geq 1 \Rightarrow \dim (\Hom_{\mathfrak
g}(V(\mu_1) \otimes V(\mu_2), V(\tau))) \geq 1.
\end{gather*}
Further, in~\cite[Theorem 1(iii)]{CFS13} the conjecture was proved for $\mathfrak{sl}_3$.
Further evidence to this conjecture is due the link between fusion products and certain PBW graded modules as provided
in~\cite{Fou14b}.

{\bf 3.6.}~Besides some numerical evidence and the partial cases stated above, the conjectures remain open.
With the result of the current paper one may tempt to replace the simple mo\-du\-les~$V(\lambda)$ in the tensor product by
minimal af\/f\/inizations of $V(\lambda)$ (see~\cite{CP96} for more details on these) and still conjecture a~surjective map
of $\mathfrak g$-modules.
Note that for classical types, ${\rm KR}(m \omega_i)$ is the minimal af\/f\/inization of $V(m \omega_i)$~\cite{Cha95, CP96}.
Some evidence to this conjecture can be found in~\cite{MY12} where the $Q$-system relations are extended to minimal
af\/f\/inizations.

\section{Proof of the main theorem}
\label{section-proof}

\begin{thm}
\label{two-factor}
Let $i\in I$ and $m_1 > m_2 > 0$.
Then there exists an embedding of $\mathfrak{g}$-modules
\begin{gather*}
{\rm KR}(m_1 \omega_i) \otimes {\rm KR}(m_2 \omega_i) \rightarrow {\rm KR}((m_1 - 1) \omega_i) \otimes {\rm KR}((m_2+1) \omega_i)
\end{gather*}
\end{thm}

The proof uses arguments of the proofs of~\cite[Theorems~1.3 and~6.1]{hcr}.
We explain below the main dif\/ferences between our present situation and the results proved in~\cite{hcr}.

\begin{proof}
For $m_1 = m_2 +1$ it is clear as the two $\mathfrak{g}$-modules are isomorphic.
Let us suppose that $m_1\geq m_2+2$.
It follows directly from the following statement that we establish.
The $U_q(\hat{\mathfrak{g}})$-module
\begin{gather*}
V = W_{m_1,q_i^{-2(m_1-m_2-1)}}^{(i)}\otimes W_{m_2,1}^{(i)}
\end{gather*}
is simple and occurs as a~composition factor in $\Rep(U_q(\hat{\mathfrak{g}}))$ of
\begin{gather*}
V' = W_{m_1-1,q_i^{-2(m_1-m_2-1)}}^{(i)}\otimes W_{m_2+1,1}^{(i)}.
\end{gather*}
For $m_1 = m_2 + 2$ this is a~direct consequence of the~$T$-system in Theorem~\ref{tsys} which gives the decomposition
of $W_{m_1-1,q_i^{-2}}^{(i)}\otimes W_{m_1-1,1}^{(i)}$ in simple modules in $\Rep(U_q(\hat{\mathfrak{g}}))$.

For a~general $m_1 > m_2 + 2$, the two modules~$V$ and $V'$ have the same highest monomial
\begin{gather*}
M = Y_{i,q_i^{2m_2}}\big(Y_{i,q_i^{2(m_2-1)}}Y_{i,q_i^{2(m_2-2)}}\cdots Y_{i,1}\big)^2\big(Y_{i,q_i^{-2}}Y_{i,q_i^{-4}}\cdots
Y_{i,q_i^{-2(m_1-m_2-1)}}\big).
\end{gather*}
Hence it suf\/f\/ices to prove that~$V$ is simple, that is $V\simeq L(M)$.
To prove it, we write a~proof as for~\cite[Theorem~6.1(2)]{hcr}:

(1) We list the dominant monomials occurring in $\chi_q(V)$, and we get as in~\cite[Lemma~5.6(2)]{hcr} the following set:
\begin{gather*}
\big\{M,MA_{i,q_i^{2m_2-1}}^{-1},MA_{i,q_i^{2m_2-1}}^{-1}A_{i,q_i^{2m_2-3}}^{-1},\dots,
MA_{i,q_i^{2m_2-1}}^{-1}A_{i,q_i^{2m_2-3}}^{-1}\cdots A_{i,aq_i}^{-1}\big\}.
\end{gather*}
All of them occur with multiplicity~$1$.

(2) We prove with the same argument as in~\cite[Section~6.1]{hcr} that these monomials do occur in $\chi_q(L(M))$:
otherwise, we would have a~monomial
\begin{gather*}
M' = MA_{i,q_i^{2m_2-1}}^{-1}A_{i,q_i^{2m_2-3}}^{-1}\cdots A_{i,aq_i^{2m_2-1-2r}}^{-1}< M
\end{gather*}
in this list such that $L(M')$ is a~simple constituent of~$V$.
But then all monomials of $\chi_q(L(M'))$ should occur in $\chi_q(L(M))$, in particular $M'A_{i,aq_i^{2m_2-1-2r}}^{-1}$,
which is not as explained in~\cite[Section~6.1]{hcr}.
\end{proof}

\pdfbookmark[1]{References}{ref}
\LastPageEnding

\end{document}